\documentclass[11pt, reqno, a4paper]{amsart}
\usepackage{soul,graphicx,color,amssymb,latexsym,amsfonts,amsmath,amsthm,enumerate,mathrsfs}
\theoremstyle{plain}

\usepackage{tikz}
\usetikzlibrary{matrix}

\newtheorem{theorem}{\bf Theorem}

\numberwithin{equation}{section} \theoremstyle{plain}
\theoremstyle{definition}

\begin{document}  




\title{Akcoglu's Dilation Theory in $L_1$-spaces}

\author{Olivia Mah}
\address{Department of Mathematics, University of San Francisco, 
San Francisco, CA 94117}
\email{ommah@usfca.edu}

\begin{abstract}
	We provide an outline for the proof of Akcoglu's dilation theory in $L_1$~\cite{Akcoglu1975a} which 
		is used by Gustafson at el.~\cite{Gustafson1997a} to investigate the embedding of a probabilistic process into a larger 			deterministic dynamical system.
\end{abstract}

\maketitle

\section{Introduction}
	
		The motivation for studying dilations in operator theory was to gain insight  
		into the properties of an operator on a smaller space by using the properties of the operator
		on a larger space.   
		
		Over the years, researchers have applied results from dilation theory to mathematical physics.  
			One example is the study conducted by Gustafson et al. in investigating the embedding 
			of a probability process into a larger deterministic dynamical 
			system~\cite{Gustafson1993, Gustafson1997a}. \footnote{R.J. Reynolds Tobacco Company 
			requested a reprint of~\cite{Gustafson1993} shortly after it was published~\cite{Gustafson2002}.}
			The physical significance of their research is related to the work on irreversible processes by Proggigne, 
			a Nobel Prize chemist~\cite{Misra1979}. 
	
		In their work~\cite{Gustafson1993, Gustafson1997a}, Gustafson et al. used two dilation theories:  one from
		Rokhlin~\cite{Cornfeld1982, Sinai1989, Rokhlin1952}, the other from Akcoglu~\cite{Akcoglu1975a, Akcoglu1975b,
			Akcoglu1976, Akcoglu1977}.  A general result was proved in~\cite{Gustafson1997a} using Akcoglu's dilation
			theorem.  
			
		With deep and interesting results, Akcoglu's dilation theorem was first proved in $L_1$ 
		and subsequently generalized to $L_p$.  The original proof in Akcoglu's dilation theorem 
			was rather complicated and simpler versions of the proof were presented~\cite{Kern1977, Nagel1982}.   
			Nonetheless, it may still be worthwhile for those who are new to the field to go through Akcoglu's proof. 
			With this in mind, the purpose of this article is to provide an outline for
			the proof in $L_1$~\cite{Akcoglu1975a}, delineating the interplay of several fundamental ideas.  				
			
		The organization of this article is as follows.  In Section 2, we present the definition of a dilation, the most general result of 
			Gustafson et al. in~\cite{Gustafson1997a} and Akcoglu's dilation theorem in $L_1$.
			In section 3, we outline the main ideas of Akcgolu's proof.  In the last section, we state
			the open problem discussed in~\cite{Gustafson1997b} and 
			present a dilation result from Rota~\cite{Rota1962}.

\section{Preliminaries}
			
	We first introduce the definition of a dilation.
			An operator $B$ on a Hilbert space $K$ is said to be a ``dilation" of an operator $A$ on a 
				Hilbert space $H$ provided $H$ is a subspace of $K$ and for all positive integers $n$,
					\[
						A^n = P B^{n} | H,
					\]	
			where P is the orthogonal projection of $K$ onto $H$ with range $H$.  (The original definition 
			introduced by Paul Halmos was a weaker form than the above as it did not 
			have the power $n$~\cite{Halmos1950}).  	
			
			A landmark theorem on dilations was proven by Sz-Nagy, which states that
		very contraction has a power dilation that is unitary~\cite{SzNagy1970}.  
		For a historical development on dilation theory, 
		see~\cite{Byrnes2001, Halmos1991, Landau1987}.  
		
		The dilation theories used to determine the types of Markov
			semigroup which could be embedded in a deterministic dynamical system 
			without approximations in~\cite{Gustafson1993, Gustafson1997a} include that of Rokhlin and 
			Akcoglu.  Roughly speaking, Sz-Nagy's dilation theory can be viewed as dilating Hilbert spaces,
			Rohlin's as dilating the underlying dynamical systems and Akcoglu's the measure spaces.  
		
		Using Akcoglu's dilation theorem, Gustafson et at. proved that an arbitrary Markov semigroup 
			satisfying certain basic properties can be dilated into a deterministic dynamical system.   
				\begin{theorem}\label{T:main_theorem}
					Every Markov semigroup $M_t$ acting on the Hilbert space $K$ 
		 				arises as a projection of a dynamical system in a larger Hilbert space  $H$.
				\end{theorem}
		
		The Markov semigroup family $\{M_t\}_{t \geq 0}$ is defined in~\cite{Gustafson1997a} as a semigroup of Markov operators,
			which are the positive integral-preserving contractions $T$ in Akcoglu's Theorem (see Theorem~\ref{T:main} below).  
			Under certain assumptions, a Markov operator is induced by transition probabilities, which in turn,
			determines a stochastic process called the Markov process (see Proposition V.4.4 in~\cite{Neveu1965}
			and~\cite{Antoniou2003}).
		
		In the next section, we state Akcoglu's dilation theorem and outline its proof.
	


\section{Akcoglu's Theorem}

	Akcoglu's dilation result (see Theorem~\ref{T:main} below) was first proved  in $L_1$~\cite{Akcoglu1975a}
		and then extended to $L_p$~\cite{Akcoglu1975b, Akcoglu1976, 
		Akcoglu1977}.  Since the ideas in the $L_p$ proof are similar to those in $L_1$, 
		for the rest of the paper, we outline the $L_1$ proof in~\cite{Akcoglu1975a}.
	
	Here is the dilation result in $L_1$~\cite{Akcoglu1975a}:
	
	\begin{theorem}\label{T:main}
		Let $(X, \mho, \mu)$ be a Borel Space and let $T$ be a positive contraction
			on $L_1(X, \mho, \mu)$.  Then there exists another Borel Space $(Y, \Im, \nu)$ 
			and a non-singular invertible transformation $\tau: Y \to Y$ so that the positive
			isometry $Q$ on $L_1 (Y, \Im, \nu)$ 
			 by $\tau$ is a dilation of $T$.
	\end{theorem}

	The setup of the proof is as follows:
	
		Let $J = [0,1]$, $\mathfrak{B}$ be the $\sigma$-algebra of the Borel subsets of $J$ and $\mu$
		is the (finite) measure on $(J, \mathfrak{B})$.   The cartesian product
		of finitely or countably many copies of $(J, \mathfrak{B})$ is denoted by: 		
			\[
				\begin{aligned}
					(J, \mathfrak{B})^n = (J^n, \mathfrak{B}^n) \quad \mbox{and} \quad 
					(J, \mathfrak{B})^{\infty} = (J^{\infty}, \mathfrak{B}^{\infty}).
				\end{aligned}
			\] 
	Also, let $(J_i, \mathfrak{B}_i)$
		be the copies of $(J, \mathfrak{B})$.
	
		Let $f \in L_1 (J^{\infty}_{-\infty}, \mathfrak{B}^{\infty}_{-\infty}, \mu^{\infty}_{-\infty})$ be a function depending
		only on $x_0$ and consider $f$ being a member of $L_1 (J, \mathfrak{B}, \mu)$ also.  Then 
		Theorem~\ref{T:main} was proved by showing that $EQ^n f = T^n f$ where

	\begin{itemize}
		\item	 $T$ is a positive contraction on $L_1 (J_0, \mathfrak{B}_0, \mu)$, 
		\item	$Q$ is the positive contraction on 
				$L_1 (J^{\infty}_{-\infty}, \mathfrak{B}^{\infty}_{-\infty}, \mu^{\infty}_{-\infty})$, and
		\item	$E$ is the conditional expectation operator with respect to $\mathfrak{B}_{0}$ i.e. 
			$E: L_1 (J^{\infty}_{-\infty}, \mathfrak{B}^{\infty}_{-\infty}, \mu^{\infty}_{-\infty}) \to L_1 (J, \mathfrak{B}_0, \mu)$.
	\end{itemize}

The two main questions that the proof has to address are:

	\begin{enumerate}
		\item	From the original measure space $X$, how do we obtain a larger space $Y$?
		\item	How do we obtain the dilation $Q$ from the larger space $Y$?
	\end{enumerate}

A short answer to the first question is to "dilate" the measure space $X$ via conditioned 
	measures (see Section~\ref{S:conditioned_measures}).  As for the second question, once 
	we obtain a larger space $Y$, we can then construct a non-singular invertible transformation $\tau$, which in turn invokes a 
	Frobenius-Perron operator serving as the dilation of $T$ (see Section~\ref{S:dilation}). We will give a more 
	detailed discussion in the next section.

\section{Outline of Proof}

\subsection{Definitions}

\subsubsection{Positive Contractions}
	A linear operator $T: L_1(S, \mathcal{S}, \vartheta) \to  L_1(S, \mathcal{S}, \vartheta)$ 
		is called a positive operator if $Tf \geq 0$ for every $f \geq 0$ 
		in $L_1(S, \mathcal{S}, \vartheta)$.  It is a contraction if its norm is less 
		than 1, that is, $|| T || \leq 1$.

\subsubsection{Transporting Measures}

	Suppose that $(S, \mathcal{S}, \vartheta)$ is a finite measure space.  Let	$\rho$ be a measurable map from
	$(S, \mathcal{S})$ to another measurable space $(S', \mathcal{S}^{'})$.  If $\rho$
	is invertible, that is, it has a measurable inverse, then 
	the set function $\varsigma$ on $\mathcal{S}^{'}$ defined by 
		\[
			\varsigma (A') = \vartheta (\rho^{-1} (A')), \quad \quad A' \in \mathcal{S'}^{'},
		\]
	is a measure in $(S', \mathcal{S}^{'})$ and is called the transport of the measure
		$\vartheta$ via $\rho$ i.e.  $\rho$ transports measure $\vartheta$ to $\varsigma$.  (If $\vartheta$
		is a probability measure and $(S', \mathcal{S}^{'}) =  (\mathbb{R}^n, \mathcal{R}^n)$, 
		then $\varsigma$ is called the distribution of $\rho$.)

		

\subsubsection{Conditioned Measures}\label{S:conditioned_measures}
	
	Let $(S, \mathcal{S})$ and $(\Xi, \mathcal{M})$ be two measurable spaces.  A family of normalized measures
		$\{\eta\}$ on $(\Xi, \mathcal{M})$ is said to be conditioned by $(S, \mathcal{S})$ 
		if this family is indexed by $s \in S$ and denoted by $\{\eta\}_{S}$.  Let $(Z, \Im) 
		=  (S, \mathcal{S}) \times (\Xi, \mathcal{M})$.
		With the conditioned family $\{\eta\}_{S}$ 
		on $(\Xi, \mathcal{M})$ and a measure
	$\vartheta$ on $(S, \mathcal{S})$, there exists a measure $\varpi = \vartheta \times \{\eta\}$ on $(Z, \Im)$ such that
		\[
			\varpi (A  \times M ) = \int_{A} \eta (M, s) \vartheta (ds),
		\]
	for each $A \in \mathcal{S}$ and $M \in \mathcal{M}$. 
	
	With $\varpi = \varphi \times \{\eta\}$, we can define the conditional expectation operator 
		$E:L_1 (Z, \Im, \varpi) \to L_1 (Z, \mathcal{S}, \varpi)$ with respect to $\mathcal{S}$ as:
			\[
				(E \, f)(s) = \int_{Y} f(s, y) \eta(dy, s), \quad  f\in L_1 (Z, \Im, \varpi), \, s \in S.
			\]

\subsubsection{Equivalence}	
	A mapping between two measurable spaces is called an equivalence if it is measurable and 
		invertible in both directions.
	
%

	

\subsection{Obtaining a Larger Space}\label{P:Larger_Space}
	
	Conditioned families of measures were investigated by both Rokhlin~\cite{Rokhlin1952} 
		and Maharam~\cite{Maharam1950}.  The Rokhlin-Maharam theorem, stated in the following 
		form in~\cite{Akcoglu1975a} and~\cite{Akcoglu1976}, provides us with tools to dilate spaces.    
		
	\begin{theorem}[Rokhlin-Maharam Theorem]\label{T:building_block}
		Let $(\Omega, \Sigma, \sigma)$ be a Borel space and let $f: \Omega \to J$ be a measurable function transporting 
			$\sigma$ to a measure $\nu$ on $(J, \mathfrak{B})$.  Then there exists 
			an isomorphism $\Phi: J^2 \to \Omega$ between $(J^2, \mathfrak{B}^2, \nu \times \{\eta \})$ and
			$(\Omega, \Sigma, \sigma)$, for some choice of the conditioned family $\{ \eta \}_J$, so that 
			$f \Phi: J^2 \to J$ is the projection of $J^2$ to its first component $J$, $\nu \times {\eta} - $ a.e.
			(If $(\Omega, \Sigma)$ is equivalent to $(J, \mathfrak{B})$, then $\Phi$ can be chosen as an equivalence.)	\end{theorem}
	
	The Rokhlin-Maharam Theorem shows how we can 
		embed $J$ (the first component) in $J^2$ and then recover it from $J^2$ via an isomorphism $\Phi$ 
		and a conditioned family $\{\eta\}_J$:
		
					\begin{tikzpicture}[node distance=2cm, auto]
  						\node (O) {$\Omega$};
 				 		\node (J) [right of=O] {$J$};
  						\node (JJ) [below of=J]{$J \times J$};
   						\draw[->] (O) to node {$f$} (J);
  						\draw[->] (JJ) to node  {$\Phi$} (O);
						\draw[->] (JJ) to node [swap] {$f \circ \Phi$} (J);
					\end{tikzpicture}

	We will briefly show how the Rokhlin-Maharam Theorem is used in the following theorem, 
		which is key to proving Akcoglu's dilation
		result in Theorem~\ref{T:main}.

		\begin{theorem}\label{T:Tf}
			Let $T$ be an integral-preserving positive contraction on $L_1 (J_0, \mathfrak{B}_0, \mu)$.  
				Then there exists a conditioned family $\{\alpha\} = \{\alpha\}_{J_0}$ on $(J_{-1},\mathfrak{B}_{-1})$ and
				an equivalence $\varphi: J_{-1} \times J_{0} \to  J_{0} \times J_{1}$ so that $\varphi$ transports
				$\{\alpha\} \times \mu$ on $J_{-1} \times J_{0}$ to $\nu \times \lambda$ on $J_{0} \times J_{-1}$
				 where $d\nu = (T1) d\mu$, so that 
					\begin{equation} \label{E:main_result}
						(Tf)(x_0) = (T1)(x_0) \int_{J_1} \, f( \varphi^{-1}_{0} (x_0, x_1)  ) \, dx_1
					\end{equation}
				for each $f \in L_1( J_0, \mathfrak{B}_0, \mu)$ and for $\mu$-a.a. $x_0 \in J_0$. 
				Here the integration is with respect to the standard Borel measure $\lambda$ on 
				$(J_{1},\mathfrak{B}_{1})$ and $\varphi^{-1}_{0} (x_0, x_1)$ denotes the $J_0$-coordinate
				of $\varphi^{-1} (x_0, x_1) \in J_{-1} \times J_{0}$, where  $(x_0, x_1) \in J_{0} \times J_{1}$. 
		\end{theorem}

	Here is how the Rokhlin-Maharam Theorem is used.  We first start with a measurable function
		$g: J_0 \times J_1 \to J$.  Then apply Theorem~\ref{T:building_block} to $g$ and obtain an equivalence
		$\varphi: J_{-1} \times J_{0} \to  J_{0} \times J_{1}$  on
		$(J_{-1}, \mathfrak{B_{-1}})$   so that $g \varphi: J_{-1} \times J_{0} \to J$ is the identification
		of $J_0$-component in $J_{-1} \times J_{0}$ via the conditioned family $\{\alpha\}_{J_0}$ on $(J_{-1},\mathfrak{B}_{-1})$
		as illustrated below:
				
		\begin{tikzpicture}[node distance=2cm, auto]
 			 \node (J0J1) {$J_{0} \times J_{1}$};
 			 \node (J) [right of=O] {$J$};
 			 \node (JNJ0) [below of=J]{$J_{-1} \times J_{0} $};
  			 \draw[->] (J0J1) to node {$g$} (J);
  			\draw[->] (JNJ0) to node  {$\varphi$} (J0J1);
  			 \draw[->] (JNJ0) to node [swap] {$g \circ \varphi$} (J);
		\end{tikzpicture}


\subsection{Integral Preserving}
	Note that the positive contraction $T$ in Theorem~\ref{T:Tf} is integral-preserving, that is,
	$\int f d\mu = \int Tf d\mu$ for all $f \in L_1$.   With the following theorem~\cite{Akcoglu1975a},
		in proving Akcoglu's dilation result, we can reduce the case of a general positive contraction 
		to that of an integral preserving positive contraction.  

	\begin{theorem}\label{T:integral_preserving}
		Every positive contraction on the $L_1$ space of a Borel space has a dilation to an operator
			of the same type which is also integral preserving.
	\end{theorem}

\subsection{Construct $\tau$}
	Recall our ultimate goal is to obtain a dilation defined on 
		$L_1 (J^{\infty}_{-\infty}, \mathfrak{B}^{\infty}_{-\infty}, \mu^{\infty}_{-\infty})$.  To do that, we 
		first need to use $\varphi$ in Theorem~\ref{T:Tf} to construct an equivalence 
		$\tau: J^{\infty}_{-\infty} \to J^{\infty}_{-\infty}$ 
		on the measurable space 
		$(J^{\infty}_{-\infty}, \mathfrak{B}^{\infty}_{-\infty})$ as follows:  Recall that $\varphi$ is defined
		as $\varphi: J_{-1} \times J_{0} \to  J_{0} \times J_{1}$.  For each $x \in J^{\infty}_{-\infty}$, 
		we define $\tau$ as:
			\[
				\begin{aligned}
					\tau_i x &= x_{i-1}, \quad \quad \mbox{if $ i \neq 0, i \neq 1$,} \\
					\tau_0 x &= \varphi_0(x_{-1}, x_0), \\
					\tau_1 x &= \varphi_1(x_{-1}, x_0).
				\end{aligned}
			\]
		Then via the conditioned measures $\{\alpha_{n}\}_{J^0_{-n+1}}$ on $(J_{-n},  \mathfrak{B}_{-n})$,
			it can be shown that $\tau$ transports the measure
			\[
				\begin{aligned}
					\mu^{\infty}_{-\infty} 
							&= \dots \times \{\alpha_{2}\} 
								\times  \{\alpha_{1}\} \times \mu \times \lambda \times \lambda \times \cdots \\
				\end{aligned}				
			\]			
		on $(J^{\infty}_{-\infty}, \mathfrak{B}^{\infty}_{-\infty}, \mu^{\infty}_{-\infty})$ to
			\[
				\begin{aligned}
					\nu^{\infty}_{-\infty} &= \dots \times \{\alpha_{2}\} 
								\times  \{\alpha_{1}\} \times \nu \times \lambda \times \lambda \times \cdots \\
				\end{aligned}				
			\]
		on the same space.   
	It follows then that 
		\begin{equation}\label{E:dv_du}
			\frac{d\nu^{\infty}_{-\infty}}{d\mu^{\infty}_{-\infty}} (\cdots, x_{-1}, x_{0},x_{1}, \cdots)
				= \frac{d\nu}{d\mu} (x_0) = (T1) (x_0),
		\end{equation}
	which is a key result to use in obtaining a dilation $Q$ as shown in the next section.
	
\subsection{Obtaining a Dilation}\label{S:dilation}

So far, we still have not shown how an integral-preserving positive contraction $T$ can invoke
	 a dilation operator $Q$.  The key to that issue
	lies with the transformation $\tau$:  as a nonsingular transformation, $\tau$ induces a positive contraction 
	called the Froebenius-Perron
	operator, which becomes the dilation operator $Q$.  
	
	Here is the definition of a Frobenius-Perron operator.

\subsubsection{Frobenius-Perron Operator}
	
	If $h: X\to X$ is a nonsingular transformation, then $h$ induces a positive contraction 
		$Q$ of $L_1 (X, \mathcal{A}, \mu)$ defined uniquely by:
		\begin{equation}\label{E:q_integral}
			\int_{h^{-1}A}  f \, d\mu = \int_{A} Qf \, d\mu
		\end{equation}
	for each $f \in L_1(X, \mathcal{A}, \mu)$ and for each $A \in \mathcal{A}$ and we call $Q$ 
	the Frobenius-Perron operator.  

		If $h$ is invertible and we let $\nu (A) = \mu (h^{-1} A)$, that is, $h$ transports measure $\mu$ to $\nu$, then
			\begin{equation}\label{E:q_operator}
				(Qf) (x) = \frac{d\nu}{d\mu} (x) f(h^{-1} x).
			\end{equation}
		(See Lasota (1994), Ch 1 \& 3 in~\cite{Lasota1994}.)

\subsubsection{$\tau$ induces Q}
	Recall that $\tau$ is defined on $(J^{\infty}_{-\infty}, \mathfrak{B}^{\infty}_{-\infty}, \mu^{\infty}_{-\infty})$.
	Since $\tau$ is invertible and non-singular
		 (because $\nu^{\infty}_{-\infty}$ is absolutely 
		continuous with respect to $\mu^{\infty}_{-\infty}$), it invokes a positive contraction $Q$
		 of  $L_1 (J^{\infty}_{-\infty}, \mathfrak{B}^{\infty}_{-\infty}, \mu^{\infty}_{-\infty})$ as in~\eqref{E:q_integral}.
		 Then by~\eqref{E:q_operator} and~\eqref{E:dv_du}, we obtain
		\begin{equation}\label{E:T_and_Q}
			\begin{aligned}
				(Qf) (\cdots, x_{-1}, x_{0},x_{1}, \cdots) 
					&= \frac{d\nu^{\infty}_{-\infty}}{d\mu^{\infty}_{-\infty}} (\cdots, x_{-1}, x_{0},x_{1}, \cdots)
							\,  f(\tau^{-1} (\cdots, x_{-1}, x_{0},x_{1}, \cdots)) \\
					&= (T1)(x_0) f(\tau^{-1} (\cdots, x_{-1}, x_{0},x_{1}, \cdots)),
			\end{aligned}
		\end{equation}
	which links the dilation $Q$ with the original operator $T$.
	Then it can be shown that $EQf = EQEf$ and eventually $EQ^n f= T^n f$.


\section{Open Question}
		
An open question was raised in~\cite{Gustafson1997b} as to whether the semigroup property of $M_t$ 
	in Theorem~\ref{T:main_theorem} can be relaxed to some wider stochastic structures 
	which are of more martingale-type or which
	permit memory effects.   

	While the open problem still remains unsolved, we state without proof an interesting result from Rota
	relating dilations and reverse martingales~\cite{Rota1962}.  
	
					
Here is Rota's dilation theorem.  

\begin{theorem}
	Let $P$ be a linear positive contraction in $L_2 (S, \sum, \mu)$ which is self-adjoint and maps
	the constant function of value 1 to 1 i.e. $P1 = 1$. Then
	\begin{enumerate}
		\item		there is a dilation of the sequence of operators $P^{2n}$ into a (reversed) martingale
				$E_n$ and
		\item		for $f$ in $L_p (S, \mathcal{F}, \mu)$, $p > 1$, $lim_{n \to \infty} P^{2n}f$
					exists almost everywhere.  
	\end{enumerate}
\end{theorem}

	Here, the dilation is in the sense that $P^{2n} = \hat{E} \circ E_n$ 
		where $E_n$ are conditional expectations of a decreasing filtration 
		and  $\hat{E}$ is a conditional expectation projecting onto $L_2 (S, \sum, \mu)$.


\end{document}